\documentclass[12pt]{amsart}
\usepackage[nohug]{diagrams}
\usepackage{latexsym}
\usepackage{amsfonts}
\usepackage{amsmath}
\usepackage{amsthm}
\usepackage{amssymb}
\usepackage{amscd}
\usepackage[mathscr]{eucal}
\usepackage[left=2.5cm,top=3cm,right=3cm,foot=2cm,headsep=1cm]{geometry}
\theoremstyle{definition}
\newtheorem{definition}[]{Definition}

\newtheorem{prop}{Proposition}
\newtheorem{thm}{Theorem}
\newtheorem{lm}{Lemma}
\newcommand{\la}{\label}
\begin{document}
\numberwithin{equation}{section}
\title[A new class of examples of group-valued moment maps]
{A new class of examples of group-valued moment maps}
%
%
    \author{Alimjon Eshmatov}
\address{Department of Mathematics,
Cornell University, Ithaca, NY 14853-4201, USA}
\email{aeshmat@math.cornell.edu}
\thanks{Research partially supported by NSF grant DMS 03-06495.}
\maketitle
\begin{abstract}
The purpose of this paper is to construct new examples of group-valued
moment maps. As the main tool for construction of such examples we use
quasi-
symplectic implosion which was introduced in \cite{HJS06}. More precisely
we show that there are certain strata of $D{\bf Sp}(n)_{\rm impl}$, the
universal
imploded space, where it is singular but whose closure is a smooth
quasi-Hamiltonian
${\bf Sp}(n) \times T$ space.
\end{abstract}
\section{Introduction}
The notion of group-valued moment map, which was introduced by Alekseev,
Malkin and Meinrenken \cite{AMM98}, is a natural generalization of classical
Hamiltonian spaces. In contrast to their classical counterpart, the moment
map takes values in a Lie group instead of the
dual of the Lie algebra. Quasi-Hamiltonian manifolds and their moment
maps share many of the features of the Hamiltonian ones, such
as reduction, cross-section and implosion.

The motivation of \cite{AMM98} for developing the theory of group-valued
moment
map came from one particularly important result.
They show that the moduli space $M(\Sigma)$ of flat connection
on a closed Riemann surface $\Sigma$ of genus $k$ is a quasi-Hamiltonian
quotient
of $G^{2k}$, which possesses natural quasi-Hamiltonian $G$-structure. 
Therefore
it is a symplectic manifold, result earlier obtained by Atiyah and Bott.
They go further
generalizing it to the case $M(\Sigma,\mathcal{C})$, the moduli space of
flat
connection on punctured Riemann surface with fixed conjugacy classes
representing homotopy classes of loops around punctures.

By imitating symplectic implosion \cite{GJS02}, J. Hurtubise, L. Jeffrey,
and R. Sjamaar
introduced the notion of group-valued implosion \cite{HJS06}. It is somewhat
similar to quasi-symplectic reduction, but instead of quotienting by the
whole
stabilizer subgroup it reduces it by its certain subgroup. While an imploded
cross-section
is almost always singular, the quasi-symplectic quotients are not. For
example,
using result of \cite{GJS02}, one can show the imploded
cross-section of $D(G)$ is singular unless the commutator subgroup of $G$ is 
a
product
of copies of ${\bf SU}(2)$. Usually they are not even orbifold unless their
universal cover
of $[G,G]$ is a product of copies of ${\bf SU}(2)$. However, like in case of
singular
quotients \cite{SL91}, using imploded cross-section theorem one
can show that imploded spaces partition into symplectic manifolds.

It was observed in \cite{HJS06} that there are certain strata of $D(G)_{\rm
impl}$
where it is singular, but whose closure is smooth. This observation lead
them to construct
new class of examples of quasi-Hamiltonian manifolds. In particular when $G$
is $A$-type
i.e. $G={\bf SU}(n)$, there is a one dimensional face of the alcove whose
corresponding stratum
has a smooth closure diffeomorphic to $S^{2n}$. Motivated by this example,
we study the
implosion for type $C$ groups, i.e. $G={\bf Sp}(n)$ unitary quaternionic
group. We show that
there is a certain stratum of $D(G)_{\rm impl}$ which has a smooth closure
diffeomorphic to
${\bf HP}^{n}$. On the other hand, it also gives new examples of
multiplicity-free
quasi-Hamiltonian spaces with non-effective $G \times T$ action.

The organization of this paper is as follows. In section 2 we recall the
definition, basic
properties and related examples of group-valued moment map. In the section 3
we review the definition and basic properties of group-valued implosion. In
this section we also give an example of ''Spinning Sphere'' constructed in
\cite{HJS06}, as motivating example of our own construction. In the section
4,
we give a construction of quasi-Hamiltonian structure on ${\bf HP}^{n}$
using an implosion.

$\mathbf {Acknowledgments.}$ It is a pleasure to express my gratitude to
my advisor, Reyer Sjamaar for suggesting
this problem and for his guidance throughout the work. I also would like to
thank Eckhard Meinrenken and Anton Alekseev for useful discussion and
valuable suggestions.

\section{Quasi-Hamiltonian Manifolds}
Let $G$ be compact, connected Lie group with Lie algebra
$\mathfrak g$. Given $G$-manifold $M$, there is induced infinitesimal Lie
algebra
action:
\begin{equation}
\xi_M(x)=\frac{d}{dt}|_{t=0}\exp(-t\xi).x \quad \textmd{for} \quad \xi \in
\mathfrak g.
\end{equation}
Recall that a Hamiltonian $G$-manifold is a symplectic $G$-manifold
$(M,\omega)$ with an equivariant map, called moment map, $\Phi:M \rightarrow
\mathfrak g^*$
satisfying relation
\begin{equation}
\iota(\xi_M)\omega=d\langle \Phi,\xi \rangle.
\end{equation}
Imitating the Hamiltonian case, \cite{AMM98} introduced the notion of
so called quasi-Hamiltonian $G$-manifolds. Recall that the
Maurer-Cartan forms $\theta_{L}, \theta_{R} \in \Omega^{1}(G,\mathfrak g)$ are defined by
$\theta_{L,g}(L(g)_{*}\xi)=\xi$ and $\theta_{R,g}(R(g)_{*}\xi)=\xi$
for $\xi \in \mathfrak g$, where $L(g)$ denotes left multiplication  
and $R(g)$ right multiplication by $g$. Let $(\cdotp,\cdotp)$ be some choice of 
inner product on $\mathfrak g$. Then
there is a closed bi-invariant three-form on $G$
$$\chi=\frac{1}{12}(\theta_{L},[\theta_{L},\theta_{L}])=
\frac{1}{12}(\theta_{R},[\theta_{R},\theta_{R}])\, .$$
\begin{definition}
A \textit{ quasi-Hamiltonian $G$-manifold } is a smooth $G$-manifold $M$
equipped with
with $G$-invariant two-form $\omega$ and $G$-equivariant map
$\Phi:M \rightarrow G$, called \textit{ group-valued moment map}, such that
the following properties hold:

(i)\ $d\omega=-\Phi^{*}\chi$

(ii)\ $ker \omega_{x}=\{\xi_{M}| \xi \in$ Ker$(Ad_{\Phi(x)}+1)\}$ for all $x
\in M$

(iii)\ $\iota(\xi_{M})\omega=
\frac{1}{2}\Phi^{*}(\theta_{L}+\theta_{R},\xi)$

\end{definition}

Basic examples of quasi-Hamiltonian manifolds are provided by conjugacy
classes and double $D(G)$. One can think of them as analogs of coadjoint
orbits and
cotangent bundle respectively.
\subsection*{Conjugacy classes}
Let $\mathcal{C}$ be a conjugacy class in $G$. Define an invariant two-form
\begin{equation}
\omega_{g}(v_{\xi},v_{\eta})= \frac{1}{2}((\eta,Ad_{g}\xi)-(\xi,Ad_{g}\eta))
\quad
\mbox{for} \quad g\in \mathcal{C},
\end{equation}
where $v_{\xi}$ and $v_{\eta}$ are fundamental vector fields induced by
conjugation action on $\mathcal{C}$. Then $(\mathcal{C},\omega)$ with the
moment
map $\Phi:\mathcal{C}\hookrightarrow G$ is a quasi-Hamiltonian $G$-space.
Moreover
$\omega$ is uniquely determined by $\Phi$.
\subsection*{Double D(G)}
It has been remarked in \cite{AMM98}, the $D(G)$ plays the same role in the
category of quasi-Hamiltonian spaces, as $T^{*}G$ does in Hamiltonian one.
As the space $D(G)$
is defined:
\begin{equation}
D(G):=G\times G.
\end{equation}
The $G\times G$ action on $D(G)$ is given by:
\begin{equation}
(g_{1},g_{2}).(u,v)= (g_{1} u g_{2}^{-1},Ad_{g_{2}}v).
\end{equation}
Define a moment map $\Phi:D(G)\longrightarrow G \times G$ \hspace{0.2cm} by
\hspace{0.2cm}
$\Phi=\Phi_{1}\times \Phi_{2}$ where
\begin{equation}
\la{mom2}
\Phi_{1}(u,v)=Ad_{u}v^{-1}, \quad \Phi_{2}(u,v)=v
\end{equation}
and two-form
\begin{equation}
\omega= -\frac{1}{2}(Ad_{v}u^{*}\theta_{L},u^{*}\theta_{L})-
\frac{1}{2}(u^{*}\theta_{L},v^{*}(\theta_{L}+\theta_{R}).
\end{equation}
The following statement was shown in \cite[Proposition 3.2]{AMM98}.
\begin{prop}
The  $(D(G), \Phi, \omega)$ is a quasi-Hamiltonian $G\times G$-manifold.
\end{prop}

Large class of examples of quasi-Hamiltonian manifolds are constructed in
\cite{AMM98} by two operations called ``Fusion'' and ``Exponentiation''.

\subsection*{Fusion.}
Unlike in Hamiltonian case, group-valued moment
maps does not behave well under restriction to subgroups or taking
products. But under slight modification one can still define these
notions in the category of quasi-Hamiltonian spaces.

\begin{thm}[Internal Fusion]
Let $(M, \omega, \Phi)$ be a quasi-Hamiltonian $G\times G\times H$-manifold,
with moment map $\Phi=\Phi_{1}\times \Phi_{2} \times \Phi_{3}: M \rightarrow
G\times G\times H$. Let $G\times H$ act via embedding
$(g,h)\rightarrow(g,g,h)$
Then $M$ equipped with the two-form $\omega +
\frac{1}{2}(\Phi^{*}_{1}\theta_{L},\Phi^{*}_{2}\theta_{R})$
and the moment map $\Phi_{1}\Phi_{2}\times \Phi_{3}: M\rightarrow G\times H$
is a quasi-Hamiltonian $G\times H$-manifold.
\end{thm}

The most important class of examples produced by this operation is fusion
product
$M=M_{1}\circledast M_{2}$, where $M_{1}$ is quasi-Hamiltonian $G$-manifold
and
$M_{2}$ quasi-Hamiltonian $G\times H$-manifold. The underlying space for $M$
is defined
as a Cartesian product $M_{1}\times M_{2}$ while the quasi-Hamiltonian 
structure
is obtained by fusing the two copies of $G$ in $G\times G \times H$.

\subsection*{Exponentiation and linearization.}
Let $(M, \omega_{0}, \Phi_{0})$
be a Hamiltonian $G$-manifold. In this section we will see how one can get
from Hamiltonian manifold a quasi-Hamiltonian one and vice-versa. First
using the inner product
on $\mathfrak{g}$  one can regard $\Phi_{0}$ as a map into $\mathfrak{g}$.
Then composing $\Phi_{0}$ with $ \exp:\mathfrak{g} \rightarrow G$ one gets
a map $\Phi: M \rightarrow G$. Then according to \cite[Proposition
3.4]{AMM98} by
slightly changing 2-form $\omega=\omega_{0}+\Phi^{*}_{0}\varpi$,
the triple $(M,\omega, \Phi)$ defines a quasi-Hamiltonian $G$-manifold.

The ``inverse functor'' called linearization, is constructed in following
way.
Let $(M, \omega, \Phi)$ be a q-Hamiltonian $G$-manifold. Moreover, assume
that
there exists an $Ad$-invariant open $U$ in $\mathfrak{g}$ such that $\exp:
U\rightarrow G$ is a diffeomorphism onto an open subset containing
$\Phi(M)$, and let
$\log: \exp U\rightarrow U$ be its inverse. Then the linearization is the
Hamiltonian
$G$-manifold $(M, \omega_{0}, \Phi_{0})$, where $\Phi_{0}=log\circ \Phi$ and
$\omega_{0}=\omega - \Phi^{\ast}_{0} \varpi$.

\subsection*{Quasi-symplectic reduction.}
One other important feature of quasi-Hamiltonian spaces is reduction.

Let $(M, \omega, \Phi)$ be a quasi-Hamiltonian $G$-manifold such that
$G$ is the product $G_{1}\times G_{2}$ where $G_{2}$ torus. Let
$\Phi=(\Phi_{1}, \Phi_{2})$ be corresponding components of moment map
$\Phi$.
We want to reduce the space with respect to the first factor. Suppose that
$g\in G_{1}$ be regular value so that $\Phi^{-1}_{1}(g)$ is a smooth
manifold. The centralizer $(G_{1})_{g}$ acts locally freely on the
submanifold $\Phi^{-1}(g)$. Then the \emph{quasi-symplectic quotient} at $g$
is defined as topological space:
\begin{equation}
M/\!/_{g}G_{1}= \Phi^{-1}_{1}(g)/(G_{1})_{g}.
\end{equation}
In good cases this quotient is a symplectic orbifold. Under above
assumptions:
\begin{thm}(\cite[Theorem 5.1]{AMM98})
The restriction of $\omega$ to $\Phi^{-1}_{1}(g)$ is
closed and $(G_{1})_{g}$-basic. The form $\omega_{g}$ on the orbifold
$M/\!/_{g}G_{1}$ induced by $\omega$ is nondegenerate. The map
$M/\!/_{g}G_{1}\rightarrow G_{2}$ induced by $\Phi_{2}$ is a
moment map for the induced $G_{2}$-action on $M/\!/_{g}G_{1}$.
\end{thm}

In case when $G_{2}$ is nonabelian $M/\!/_{g}G_{1}$ is not symplectic, but
a quasi-Hamiltonian $G_{2}$-orbifold.

In singular case the quotient stratifies into symplectic manifolds according
to orbit type. Let $g$ be a arbitrary element of $G_{1}$. For each subgroup
$H$
define a $(G_{1})_{g}$-invariant submanifold $M_{(H)}$ consisting of all
points
such that the stabilizer $(G_{1})_{g}\cap(G_{1})_{x}$ is conjugate to $H$.
Put $Z=\Phi^{-1}_{1}(g)$ and $Z_{(H)}=Z\cap M_{(H)}$. Let ${Z_{i}}$ be the
collection of connected components of $Z_{(H)}$, where $H$ ranges over
conjugacy
classes of $(G_{1})_{g}$. Then we have decomposition:
\begin{equation}
\la{dec}
M/\!/_{g}G_{1}=\coprod_{i\in I}Z_{i}/(G_{1})_{g}.
\end{equation}
\begin{thm}(\cite{HJS06}) The decomposition \eqref{dec} is a locally
normally trivial stratification of $M/\!/_{g}G_{1}$ into
symplectic submanifolds. Moreover, the stratification is
$G_{2}$-invariant and the continuous map $\bar{\Phi}_{2}:
M/\!/_{g}G_{1} \rightarrow G_{2}$ induced by $\Phi_{2}$ restricts
to a moment map for the $G_{2}$ action on each stratum.
\end{thm}

\section{Group-valued imploded cross-section.}
Let $G$ be a simply connected
compact Lie group with maximal torus $T$. Recall that the symplectic
implosion is an ``abelianization functor'', which transforms a Hamiltonian
$G$-manifold into a Hamiltonian $T$-space preserving some of properties of
the
manifold, but in the expense of producing singularities \cite{GJS02}.
However,
the singularities are not arbitrary in the sense that it ``stratifies'' into
symplectic submanifolds in such a way that $T$ action preserves the
stratification.

Now let $(M,\omega,\Phi)$ be a quasi-Hamiltonian $G$-space. In
\cite{AMM98}, it was shown that like in Hamiltonian case one can
prove a convexity theorem. But the moment map image instead of a
Weyl chamber, one have to consider in an alcove. Here the
assumption of simply connectedness of the group is crucial, since
otherwise the description of a space of conjugacy classes is quite
complicated.

Let $\mathcal{C} ^{\vee}$ be the chamber in $\mathfrak{t}$ dual to
$\mathcal{C}$
and let $\mathcal{A}$ be the unique (open) alcove contained in $\mathcal{C}
^{\vee}$
such that $0 \in \bar{\mathcal{A}}$. Using the exponential map one can
identify
$\bar{\mathcal{A}}$ with space of conjugacy
classes
$T/W\cong G/AdG$, where $W$ is the corresponding Weyl group. Let denote by
$G_{g}$ the centralizer of $g$ in $G$. For points $m_{1}, m_{2}\in \Phi^{-1}
(\exp\bar{\mathcal{A}})$ define $m_{1}\sim m_{2}$ if $m_{2}=gm_{1}$ for some
$g \in [G_{\Phi(m_{1})},G_{\Phi(m_{1})}]$. One can check that $\sim$ is
indeed
equivalence relation.
\begin{definition}
The imploded cross-section of $M$ is the quotient space $M_{\rm impl}=
\Phi^{-1}(\exp\bar{\mathcal{A}})/\sim$, equipped with the quotient
topology. The imploded moment map $\Phi_{\rm impl}$ is the continuous
map $M_{\rm impl}\rightarrow T$ induced by $\Phi$.
\end{definition}
$M_{\rm impl}$ has many nice properties that smooth manifolds posses.
It is Hausdorff, locally compact and second countable. The action of
$T$ preserves $\Phi^{-1}(\exp\bar{\mathcal{A}})$ and descends to a
continuous action on $M_{\rm impl}$.

We have decomposition of $M_{\rm impl}$ into orbit spaces:
\begin{equation}
\la{decomp}
M_{\rm
impl}=\coprod_{\sigma\leq\mathcal{A}}\Phi^{-1}(\exp\sigma)/[G_{\sigma},
G_{\sigma}],
\end{equation}
where $\sigma$ ranges over the faces of alcove $\mathcal{A}$ and
$K_{\sigma}$
is the centralizer of $\exp\sigma$. Let us denote piece
$\Phi^{-1}(exp\sigma)/[G_{\sigma},G_{\sigma}]$ by $X_{\sigma}$. Using
quasi-Hamiltonian cross-section [HJS06, Theorem 3.4] each $X_{\sigma}$
stratifies into symplectic manifolds. Now let $\{X_i|i\in I\}$ be the 
collection of the all strata of the all pieces $X_{\sigma}$ where $\sigma$ 
ranges over the faces of alcove. Then imploded cross-section $M_{\rm
impl}$ the disjoint union:
\begin{equation}
\la{dec1}
M_{\rm impl}=\coprod_{i\in I}X_{i}
\end{equation}
such that each piece $X_{i}$ is symplectic manifold:
\begin{thm}(\cite{HJS06})
The decomposition \eqref{dec1} of the imploded cross-section is a locally
finite partition into locally closed subspaces, each of which is a
symplectic manifold. There is a unique open stratum, which is dense,
in $M_{\rm impl}$ and symplectomorphic to the principal cross section of
$M$.
The action of the maximal torus $T$ on $M_{\rm impl}$ preserves the
decomposition and the imploded moment map $\Phi_{\rm impl}:M_{\rm impl}
\rightarrow T$ restricts to a moment map for the $T$-action on each
stratum.
\end{thm}
Therefore we call $M_{\rm impl}$ a \textit{stratified quasi-Hamiltonian
$T$-space}.

\subsection*{Imploded cross-section of the double.}
In the example of q-Hamiltonian manifolds we have seen that $D(G):=G\times
G$
possesses quasi-Hamiltonian $G\times G$-structure.

Let $M$ be an arbitrary quasi-Hamiltonian $G$-space. Fusing it with $D(G)$
one obtains a q-Hamiltonian $G\times G$-manifold $M\circledast D(G)$.
Now define $j: M\rightarrow M\circledast D(G)$ by $j(m)= (m,1,\Phi(m))$.
Then one of the main results of \cite{HJS06} states:
\begin{thm}[universality of imploded double]
Let $M$ be a quasi-Hamiltonian $G$-manifold. The map $j$ induces a
homeomorphism
\begin{equation}
j_{\rm impl}:M_{\rm impl}\xrightarrow{\cong}(M\circledast D(G)_{\rm
impl})/\!/G
\end{equation}
which maps strata to strata and whose restriction to each stratum is an
isomorphism of quasi-Hamiltonian $T$-manifolds.
\end{thm}

\subsection*{A smoothness criterion and quasi-Hamiltonian structure on
$S^{2n}$}
We have seen that in previous section that in order to construct implosion
of
a given manifold, it suffices to know the implosion of double of
corresponding
Lie group. The implosion of double is singular space, however the
singularities
on certain strata are removable. In order to show that one has to use the
explicit
correspondence between $D(G)$ and $T^{*}G$.

Identify $\mathfrak{g}$ with $\mathfrak{g}^{*}$ using bi-invariant inner
product on $\mathfrak{g}$. Trivializing $T^{*}G$ in a left-invariant manner,
define $G\times G$-equivariant  map $\mathcal{H}=id\times \exp: T^{*}G
\rightarrow D(G)$. Let $(T^{*}G, \omega_{0}, \Psi_{0})$ be a Hamiltonian
$G\times G$ manifold, where $\omega_{0}$ is the canonical symplectic form
on cotangent bundle and a moment map
$\Psi_{0}(g,\lambda)=(-Ad_{g}\lambda,\lambda)$.
Let $O$ be the set of all $\xi \in \mathfrak{t}$ with $(2\pi
i)^{-1}\alpha(\xi)<1$
for all positive roots $\alpha$ and $U=(AdG)O$.
\begin{lm}(\cite{HJS06})
The triple $(T^{*}G, \mathcal{H}^{*}\omega, \mathcal{H}^{*}\Psi)$ is the
exponentiation of $(T^{*}G, \omega_{0}, \Psi_{0})$. In particular, $G\times
U$
is a quasi-Hamiltonian $G\times G$-manifold.
\end{lm}
Now using local diffeomorphism given by $\mathcal{H}$ and that of result of
\cite{GJS02} we have:
\begin{thm}[Smoothness criterion]
Let $\sigma$ be a face of $\mathcal{A}$ satisfying $[G_{\sigma},G_{\sigma}]
\cong \mathbf{SU}(2)^{k}$ (resp. $[\mathfrak{g}_{\sigma},
\mathfrak{g}_{\sigma}]
\cong \mathfrak{su}(2)^{k}$) for some $k\geq0$ and possessing a vertex $\xi$
such that $\exp \xi$ is central. Then $D(G)_{\rm impl}$ is a smooth
quasi-Hamiltonian
$G\times T$-manifold(resp. orbifold) in a neighborhood of the stratum
corresponding
to $\sigma$.
\end{thm}
The partial converse of this result is also true. Suppose that $\sigma$
contains
a vertex $\xi$ such that $\exp \xi$ is central and $D(G)_{\rm impl}$ is 
smooth
in a
neighborhood of the corresponding stratum. Then
$[G_{\sigma},G_{\sigma}]\cong
\mathbf{SU}(2)^{k}$. On the other hand, there are certain strata where
$D(G)_{\rm impl}$
is singular, but their closure is a smooth quasi-Hamiltonian manifold.

Let $G$ be $\mathbf{SU}(n)$. Consider an edge $\sigma_{01}$ of an alcove
with
centralizer $G_{01}=\mathbf{S}(\mathbf{U}(1)\times \mathbf{U}(n-1))$. By the
argument
above we know that for $n>3$ the corresponding stratum $X_{01}$ in $X$
consists of
genuine singularities. Nevertheless the following result asserts that it is
smooth
quasi-Hamiltonian manifold and in fact $S^{2n}$.
\begin{thm}(\cite{HJS06})
The closure of the stratum $X_{01}$ of $X=D\mathbf{SU}(n)_{\rm impl}$ is a
smooth
quasi-Hamiltonian $\mathbf{U}(n)$-manifold diffeomorphic to $S^{2n}$.
Furthermore
antipodal map of $S^{2n}$ corresponds to involution of $X_{01}$ obtained by
lifting
symmetry of the alcove $\mathcal{A}$ that reverses the edge $\sigma_{01}$.
\end{thm}

\section{Imploded cross-section of $\mathbf{Sp}(n)$}
In the last section, we saw a construction using an imploded cross-section
yielding
an example of sphere with quasi-Hamiltonian structure. In this section using
construction
of somewhat similar nature, we will show that $\mathbf{HP}^{n}$ have
quasi-Hamiltonian
structure as well.

Let $G=\mathbf{Sp}(n)$, the group of unitary $n\times n$ matrices over the
quaternions,
with maximal torus $T=\{\mathtt{diag}(e^{2\pi i x_{1}},...,e^{2\pi i
x_{n}})\}$.
Identify $\mathfrak{t}$ with $\mathbf{R}^{n}$ via the map $x\mapsto 2\pi
i\,\mathtt{diag}
(x_{1},...,x_{n})$. Then the simple roots will have form:
\begin{equation}
(2\pi i)^{-1}\alpha_{k}(x)=x_{k}-x_{k+1} \quad \textmd{for} \quad
k=1,...,n-1 \quad \textmd{and} \quad (2\pi i)^{-1}\alpha_{n}(x)=-2x_{1}
\end{equation}
with minimal root $(2\pi i)^{-1}\alpha_{n}(x)=2x_{n}$. The
corresponding alcove is the $n$-simplex $0<x_{n}<...<x_{1}<1/2$.
We will abuse our notation and denote $\sigma_{01}$ the edge of
the simplex with vertices $\sigma_{0}$, $\sigma_{1}$ corresponding
to $I$ and $\mathtt{diag}(-1,1,...,1)$. Under exponential map this
edge corresponds to torus elements of the form
$\mathtt{diag}(t,1,...,1)$ with centralizer
$G_{01}=\textbf{U}(1)\times\textbf{Sp}(n-1)$. The centralizers of
the vertices are $G_{0}=\mathbf{Sp}(n)$ and $G_{1}=\mathbf{Sp}(1)
\times \mathbf{Sp}(n-1)$ respectively. By \eqref{decomp}, the
corresponding strata are given by:
\begin{equation}
X_{\sigma}=G/[G_{\sigma},G_{\sigma}]\times \exp \sigma.
\end{equation}
Therefore
\begin{equation}
X_{0}=\{\overline{I}\}\times \{I\}\cong \{\mathtt{pt}\},
\end{equation}
\begin{equation}
X_{01}=\mathbf{Sp}(n)/\mathbf{Sp}(n-1)\times \{\mathtt{diag}(t,1,...,1)|t\in (0,1)\}
\cong S^{4n-1}\times(0,1),  
\end{equation}
\begin{equation}
X_{1}=\mathbf{Sp}(n)/(\mathbf{Sp}(n-1)\times\mathbf{Sp}(1))\times
\mathtt{diag}(-1,1,...,1)\cong
\mathbf{HP}^{n-1}.
\end{equation}
where overline denotes a coset corresponding to an element.
Consider the closure of the stratum corresponding to $\sigma_{01}$,
$\bar{X}_{01}=
X_{0}\bigcup X_{01} \bigcup X_{1}$. Notice that $X_{0}\bigcup X_{01}\cong
\mathbf{H}^{n}$
and $X_{1}\cong \mathbf{HP}^{n-1}$, and one would expect that
$\bar{X}_{01}\cong \mathbf{HP}^{n}$.
We will prove this by directly constructing a homeomorphism from
$\bar{X}_{01}$ to $\mathbf{HP}^{n}$.
Define a map:
\begin{equation}
\mathcal{G}: \bar{X}_{01} \rightarrow \mathbf{HP}^{n} \quad , \quad
( \bar{g}, x)\mapsto \Bigr[\sqrt{\frac{2x_{1}}{1-2x_{1}}},g.v \Bigr],
\end{equation}
where $v=(1,0,...,0) \in \mathbf{H}^{n}$. One can easily check that it is
well-defined,
i.e., does not depend on the equivalence class of $g$ in
$\mathbf{Sp}(n)/\mathbf{Sp}(n-1)$.
$\mathcal{G}$ is a continuous, bijective on $X_{01}$(or $0<x_{1}<1$) and
continuously extends to
$\bar{X}_{01}$. Indeed on $X_{0}$ (or $x_{1}=0$) $\mathcal{G}( \bar{g},
x)=[1,0,...,0]$ and
on $X_{1}$ (or $x_{1}=1$) $\mathcal{G}( \bar{g}, x)=[0,g.v]$.

In \cite{HJS06} it was shown that imploded space is Hausdorff, locally
compact and second
countable. Therefore $\mathcal{G}$ is a homeomorphism. Define a smooth
structure on
$\bar{X}_{01}$ by pulling back a smooth structure on $\mathbf{HP}^{n}$ via
homeomorphism
$\mathcal{G}$. Then the inverse map $\mathcal{F}:\mathbf{HP}^{n} \rightarrow
\bar{X}_{01}$
is smooth and defined as:
\begin{equation}
([Z_{1},...,Z_{n+1}]) \mapsto\overline{\left(\begin{array}{c} f_{p,n-q}(Z)
\\\end{array}\right)}
\times \mathtt{diag}(e^{2\lambda \pi i},1,...,1),
\end{equation}
where $\lambda$ and $f_{p,n-q}(Z)$ are given by:
\begin{equation}
\lambda =\frac{\sum_{l=2}^{n+1}|Z_{l}|^{2}}{\sum_{l=1}^{n+1}|Z_{l}|^{2}}
\end{equation}
\begin{equation}
f_{p,n-q}=
\left \{
\begin{array}{ccccc}
   \frac{Z_{p+1}\bar{Z}_{1}}{|Z_{1}|\sqrt{\sum_{l=2}^{n+1}|Z_{l}|^{2}}}
   &\big \vert &  \mbox{if} &  q=n-1
\\*[2ex]

\frac{|Z_{q+3}|Z_{p+1}}{\sqrt{\sum_{l=2}^{q+2}|Z_{l}|^{2}}\sqrt{\sum_{l=2}^{q+3}|Z_{l}|^{2}}}
      &\big \vert & \mbox{if} & p-q<2   \mbox{ and} & q\leq n-2
\\*[2ex]

\frac{\sqrt{\sum_{l=2}^{p}|Z_{l}|^{2}}Z_{p+1}}{\sqrt{\sum_{l=2}^{p+1}|Z_{l}|^{2}}|Z_{p+1}|}
      &\big \vert & \mbox{if} & p-q=2   \mbox{ and} & q\leq n-2
\\*[2ex]
    0
  &\big \vert &  \mbox{if} & \mbox{ otherwise} &
\end{array}
\right. \ ,
\end{equation}
There is an action of $\mathbf{Sp}(n) \times T$ on $\mathbf{HP}^{n}$ :
\begin{equation}
(g,t)[Z_{1},...,Z_{n+1}]=[Z_{1}t_{1}, k.(Z_{2},...,Z_{n+1})] \ ,
\end{equation}
where we regard $\mathbf{H}^{n+1}$ as a right $\mathbf{H}$-module.
Then one can easily show
\begin{lm}
The map $\mathcal{G}$ is a $\mathbf{Sp}(n) \times T$-equivariant.
\end{lm}
First, let us recall the quasi-Hamiltonian structure on stratum
$X_{\sigma}$.
Since the moment map $\Phi_{2}$ defined as in \eqref{mom2} is transversal to
all faces
of the alcove, using quasi-symplectic reduction with quasi-Hamiltonian
cross-section
theorem one can show:
\begin{lm}(\cite{HJS06})
For every $\sigma \leq \mathcal{A}$ the subspace $X_{\sigma}=G/[G_{\sigma},
G_{\sigma}] \times \exp \sigma$ of $D(G)_{\rm impl}$ is a quasi-Hamiltonian
$G \times T$-manifold. The moment map $X_{\sigma}\rightarrow G\times T$ is
the
restriction to $X_{\sigma}$ of the continuous map $\Phi_{\rm
impl}\rightarrow
G\times T$ induced by $\Phi: D(G) \rightarrow G\times G$.
\end{lm}
Next we compute the corresponding 2-form $\omega_{\sigma}$ on $X_{\sigma}$.
Let
$(g,x)$ be an arbitrary point on $X_{\sigma}$. A tangent vector at $(g,x)$
is of
the form $((L_{g})_{*}\xi, (L_{x})_{*}\eta)$ where $\xi \in \mathfrak{g}$
and
$\eta \in \zeta+ \mathfrak{z}(\mathfrak{g}_{\sigma})$ for some $\zeta\in
\mathfrak{z}(\mathfrak{g}_{\sigma})^{\bot}=[\mathfrak{g}_{\sigma},\mathfrak{g}_{\sigma}]$
(Lemma A.3, \cite{HJS06}). Then a simple calculation yields:
\begin{eqnarray}
\la{fm1}
(\omega_{\sigma})_{(g,x)}(((L_{g})_{*}\xi_{1},
(L_{x})_{*}\eta_{1}),((L_{g})_{*}\xi_{2},
(L_{x})_{*}\eta_{2}))  \\ \nonumber
=-\frac{1}{2}((Ad_{x}-Ad_{x^{-1}})\xi_{1},\xi_{2})-(\xi_{1},\eta_{2})+(\xi_{2},\eta_{1}).
\end{eqnarray}
One can check that it does not depend on equivalence class of $\xi_{i}$ in
$\mathfrak{g}/
[\mathfrak{g}_{\sigma},\mathfrak{g}_{\sigma}]$. Consider 2-form
$\omega_{01}$ on an open stratum
$X_{01}$. In what follows we compute the pull back of this 2-form on via
$\mathcal{F}$ and
show that it extends smoothly to all of $\mathbf{HP}^{n}$. Since
$\omega_{01}$ is
$\mathbf{Sp}(n) \times T$-invariant, it suffices to consider vectors of the
form
$z_{0}=[t,1,0,...,0]$, where
\begin{equation}
t=\frac{|Z_{1}|} {\sqrt{\sum_{l=2}^{n+1}|Z_{l}|^{2}}}. \, .\nonumber
\end{equation}
The tangent space at $z_{0}$ will be:
\begin{equation}
\la{tan}
T_{z_{0}}\mathbf{HP}^{n}=\{(w_{1},...,w_{n+1})\in\mathbf{H}^{n+1}|t
w_{1}+w_{2}=0\},
\end{equation}
where $w_{l}=w_{l1}+w_{l2}i+w_{l3}j+w_{l4}k$.
Let us first find corresponding tangent vectors at $\mathcal{F}(z_{0})$, or
more precisely
corresponding pull-backs $\xi_{i}, \eta_{i}$ to elements of Lie algebra as
in \eqref{fm1}.
Let $v$ and $w$ be tangent vectors of form \eqref{tan}. Note, since
\begin{equation}
\mathcal{F}(z_{0})=(\bar{I}, \mathtt{diag}(\exp(\lambda\pi
i),1,...,1))=:(g,x),
\nonumber
\end{equation}
the first component of the image is already an element of Lie algebra, while
the second one
has to be translated by an appropriate element of Lie group (that is $x$).
Denote by
$(A^{v}_{p,q})$ and $(B^{v}_{p,q})$ ($(A^{w}_{p,q})$ and $(B^{w}_{p,q})$ )
the matrix
representation of $\xi_{1}$ and $\eta_{1}$ (correspondingly $\xi_{2}$ and
$\eta_{2}$).
Then substituting these to the first term of \eqref{fm1} and using above
expression for
$\mathcal{F}(z_{0})$ we have:
\begin{eqnarray}
\la{fm2}
\nonumber
&&((Ad_{x}-Ad_{x^{-1}})\xi_{1},\xi_{2})
= \textrm{Re} \Big(\big[\exp(\lambda\pi i)A_{11}^{v}\exp(-\lambda\pi
i)-\\*[1ex]
   && \exp(-\lambda\pi i)A_{11}^{v}\exp(\lambda\pi i)\big]\bar{A}_{11}^{w}
   - \big[\exp(\lambda \pi i)-\exp(-\lambda \pi i)\big]
\sum_{p=2}^{n}A_{1p}^{v}A_{p1}^{w}+\\*[1ex]
&& \sum_{p=2}^{n}A_{1p}^{v}\big[\exp(-\lambda \pi i)-\exp(\lambda \pi
i)\big]
\bar{A}_{p1}^{w}\Big) \, ,
\nonumber
\end{eqnarray}
where the inner product is given by $(A,B)={\rm Re}(tr(A\bar{B}^{t}))$ and
\begin{equation}\la{dd1}
A_{11}^{v}= -(t+t^{-1})(v_{12} i + v_{13} j + v_{14} k) \, ,
\end{equation}
and by skew-symmetry:
\begin{equation}\la{dd2}
A_{p1}^{v}= - A_{1p}^{v}= v_{(p+1)1} + v_{(p+1)2} i + v_{(p+1)3} j +
v_{(p+1)4} k \, .
\end{equation}
There are similar relations to \eqref{dd1} and \eqref{dd2} if we replace $v$
by $w$.
Thus we can rewrite  \eqref{fm2} in the following form
\begin{eqnarray}
((Ad_{x}-Ad_{x^{-1}})\xi_{1},\xi_{2})=
2 \sin(2\pi \lambda) (t + t^{-1})(v_{13}w_{14}-w_{13}v_{14})- \\*[1ex]
\nonumber
4 \sin(\pi \lambda)
\sum_{p=3}^{n+1}(v_{p1}w_{p2}-w_{p1}v_{p2}-v_{p3}w_{p4}+w_{p3}v_{p4}).
\end{eqnarray}
Hence, corresponding two-form will be:
\begin{equation}
\la{eq3.16}
2\sin(2\pi
\lambda)(t+t^{-1})dx_{13}dx_{14}-4\sin(\pi\lambda)\sum_{p=3}^{n+1}(dx_{p1}dx_{p2}-
dx_{p3}dx_{p4}) ,
\end{equation}
where $x$'s are just real coordinates for $Z$'s, such that
$Z_{l}=x_{l1}+x_{l2}+x_{l3}+x_{l4}$.
For the remaining part of \eqref{fm1} we have:
\begin{equation}
\la{eq3.17}
-(\xi_{1},\eta_{2})+(\xi_{2},\eta_{1}) = {\rm
Re}(-A_{11}^{v}\bar{B}_{11}^{w}+A_{11}^{w}\bar{B}_{11}^{v}) ,
\end{equation}
where
\begin{equation}
\la{eq3.18}
B_{11}^{v}=-\frac{2\pi i t}{t^{2}+1}v_{11} ,
\end{equation}
and corresponding two-form will be:
\begin{equation}
\la{eq3.19}
2\pi dx_{11}dx_{12}.
\end{equation}
Combining \eqref{eq3.16} with \eqref{eq3.19} yields:
\begin{eqnarray}
\la{pb}
& \mathcal{F}^{*}\omega_{01}=2\pi
dx_{11}dx_{12}-2\sin(2\pi\lambda)(t+t^{-1})
dx_{13}dx_{14}+ \\*[1ex] \nonumber
& 4\sin(\pi\lambda)\sum_{p=3}^{n+1}(dx_{p1}dx_{p2} - dx_{p3}dx_{p4}) .
\end{eqnarray}
It is a smooth two-form defined on open dense subset of $\mathbf{HP}^{n}$.
Moreover we can show
\begin{lm}
\la{lem4}
The two-form $\mathcal{F}^{*}\omega_{01}$ extends smoothly on
all of $\mathbf{HP}^{n}$.
\end{lm}
\begin{proof}
It suffices to check two critical cases $Z_{1}=0$, a line at infinity,
and $[1,0,...,0]$, a point at infinity. As $|Z_{1}|$ approaches to $0$,
$\lambda$ tends to $1$ and therefore the third expression on the right
hand side of \eqref{pb} vanishes. Now since $\lambda = 1-t^{2}$,
we have $\, t \rightarrow 0\, $ and hence
$$   2\sin(2\pi\lambda)(t+t^{-1})\longrightarrow -4\pi  \, .$$
So in the neighborhood of $Z_{1}=0$ the two-form
$\mathcal{F}^{*}\omega_{01}$
can be written as
$$2\pi dx_{11}dx_{12} + 4\pi dx_{13}dx_{14} .$$
In the similar fashion one can show that in the neighborhood of
$[1,0,...,0]$
it is given by:
$$2\pi dx_{11}dx_{12} - 4\pi dx_{13}dx_{14} .$$
This finishes the proof of this lemma.
\end{proof}

Notice that the obtained 2-form
is given in dehomogenized coordinates:
\begin{equation}
\la{bb}
\Bigg[\frac{|Z_{1}|}{\sqrt{\sum_{l=2}^{n+1}|Z_{l}|^{2}}},
\frac{Z_{2}\bar{Z}_{1}}{|Z_{1}|\sqrt{\sum_{l=2}^{n+1}|Z_{l}|^{2}}},
.....,\frac{Z_{n+1}\bar{Z}_{1}}{|Z_{1}|\sqrt{\sum_{l=2}^{n+1}|Z_{l}|^{2}}}\Bigg].
\end{equation}
For $Z=x_{1}+x_{2}i+x_{3}j+x_{4}k$ define ${\rm Im}_{i}(Z)=x_{2}$, then we
have:
\begin{equation}
\la{cc}
{\rm Im}_{i}(d\bar{Z_{p}}dZ_{p})= 2(dx_{p1}dx_{p2}-dx_{p3}dx_{p4}) .
\end{equation}
Now using \eqref{bb} and \eqref{cc} in homogenous coordinates the first two
terms vanishes,
our 2-form will take form:
\begin{eqnarray}
\la{form}
\nonumber
&&
4\sin(\lambda\pi)\Bigg(|Z_{1}|^{2}\sum_{l=2}^{n+1}|Z_{l}|^{2}\Bigg)^{-1}\Bigg[
\sum_{p=3}^{n+1}|Z_{p}|^{2}{\rm Im}_{i}(dZ_{1}d\bar{Z}_{1})-
\\*[1ex]
&& {\rm Im}_{i}(Z_{1}d\bar{Z}_{p}dZ_{p}\bar{Z}_{1})+
\Bigg(\sum_{p=3}^{n+1}|Z_{p}|^{2}{\rm Im}_{i}((Z_{1}d\bar{Z}_{1}))-
  {\rm Im}_{i}(Z_{1}\bar{Z}_{p}dZ_{p}\bar{Z}_{1})\Bigg) \times \\*[1ex]
  \nonumber
&& \Bigg(\frac{Z_{1}d\bar{Z}_{1}+dZ_{1}\bar{Z}_{1}}{|Z_{1}|^{2}}+
  \frac{\sum_{l=2}^{n+1}(Z_{l}d\bar{Z}_{l}+dZ_{l}\bar{Z}_{l})}
  {\sum_{l=2}^{n+1}|Z_{l}|^{2}}\Bigg)\Bigg] .
\end{eqnarray}

Last thing we need to show that there is well-defined smooth moment map.
Define
a map $\Phi: \mathbf{HP}^{n} \rightarrow \mathbf{Sp}(n) \times T$ such that
the following diagram commutes:
\begin{equation}
\la{maps}
\begin{diagram}[small, tight]
        &           &      X_{\sigma} &  &
        \\
        &\ldInto^{\mathcal{G}|_{X_{\sigma}}} &   &\rdTo^{\Phi_{\sigma}}   \\
\mathbf{HP}^{n} &                 &\rDashto^{\Phi}   & & \mathbf{Sp}(n)
\times T \, ,\\
\end{diagram}
\end{equation}
for each face $\sigma$ in the closure. Then it has to be of the form:
\begin{equation}
\la{mom}
[Z_{1},...,Z_{n+1}] \mapsto (A B^{-1} A^{-1}, B)
\end{equation}
where $A= (A_{pq})$ and $B= (B_{pq})$ are matrices:
$$
A_{p1}=\frac{Z_{p+1}\bar{Z}_{1}}{|Z_{1}|\sqrt{\sum_{l=2}^{n+1}|Z_{l}|^{2}}}
$$
\begin{equation}
B=\mathtt{diag}(\exp(\lambda\pi i),1,...,1)
\end{equation}
Notice that we only defined the first column of the $A_{pq}$ since it is
uniquely
determined by its first column. Evidently, $\Phi$ is uniquely determined and
$\mathbf{Sp}(n)
\times T$-equivariant. We have to show that it is smooth. From the
construction one can see
that $B_{pq}$ are smooth. As for the first component of $\Phi$, using the
fact
$A \in \mathbf{Sp}(n)$ we have:
$$
A B^{-1} A^{-1} = \textmd{Id}_{n} + C ,
$$
where $C=(C_{pq})$:
$$
C_{pq}= A_{p1} \bar{B}_{11} \bar{A}_{q1} - A_{p1} \bar{A}_{q1} ,
$$
or to be more precisely:
\begin{equation}
C_{pq}=\Bigg(|Z_{1}|^{2}\sum_{l=2}^{n+1}|Z_{l}|^{2}\Bigg)^{-1} Z_{p+1}
\Big[\bar{Z}_{1} \exp(\pi i \lambda) Z_{1}-|Z_{1}|^{2}\Big]\bar{Z}_{q+1} .
\end{equation}
We can easily see that it is smooth for $Z_{1}\neq 0$ and
$\sum_{l=2}^{n+1}|Z_{l}|^{2} \neq 0$. Otherwise using almost the
same argument as in Lemma \ref{lem4} we can show it is smooth in these two
cases
as well. Now summarizing these facts we have:
\begin{thm}
The closure of the stratum $X_{01}$ of $X=D\mathbf{Sp}(n)_{\rm impl}$ is a
smooth quasi-Hamiltonian
$\mathbf{Sp}(n) \times T$-manifold diffeomorphic to $n$-dimensional
quaternionic projective
space with 2-form and moment map determined by \eqref{form} and \eqref{mom}
correspondingly.
\end{thm}

\begin{thebibliography}{1}
%
\bibitem[AMM98]{AMM98} A. Alekseev, A. Malkin, and E. Meinrenken, \textit{
Lie group valued
moment map}, J. Differential Geom. \textbf{48} (1998), 445-495. MR
99k:58062.
%
\bibitem[AMW02]{AMW02}
A. Alekseev, E. Meinrenken and C. and Woodward, \textit{ Duistermaat-Heckman
measure
and moduli spaces of flat bundles over surfaces}, Geom. Funct. Anal.
\textbf{12} (2002), 1-31. MR 2003d:53151.
%
\bibitem[GJS02]{GJS02} V. Guillemin, L. Jeffrey, and R. Sjamaar, \textit{
Symplectic implosion},
Transform. Groups \textbf{7} (2002), 155-184. MR 2003b:53090.
%
\bibitem[HJ00]{HJ00} J. Hurtubise and L. Jeffrey, \textit{ Representation
with weighted
frames and framed parabolic bundles}, Canad.J.Math. \textbf{52} (2000),
1235-1268.
%
\bibitem[HJS06]{HJS06} J. Hurtubise, L. Jeffrey, and R. Sjamaar, \textit{
Group-valued
implosion and parabolic structures}, Amer.J.Math. \textbf{128} (2006),
167-214. 2007d:53141.
%
\bibitem[SL91]{SL91} R. Sjamaar and E. Lerman, \textit{ Stratified
symplectic spaces and
reduction}, Ann. of Math. (2) \textbf{134} (1991), 375-422. MR 92g:58036.
%
\end{thebibliography}
\end{document}